\documentclass[12pt,twoside,dvips]{amsart}
\usepackage{math,amsfonts,pb-diagram,lamsarrow,pb-lams,epsfig,fancyheadings}

\newcommand\modop{\operatorname{mod}}
\newcommand\gr{{\mathfrak G}}
\newcommand\grt{{\widetilde{\mathfrak G}}}

\makeatletter
\global\renewcommand\dggeometry{\dg@XGRID=\thr@@
  \dg@YGRID=\@ne
  \unitlength=0.015pt\relax}
\makeatother

\begin{document}
\title{Lie Methods in Growth of Groups and Groups of Finite Width}
\author{Laurent Bartholdi \and Rostislav I. Grigorchuk}
\date{}
\thanks{The authors express their thanks to the Swiss National Science Foundation; the second author thanks the Russian Fund for Fundamental Research, research grant 01-00974 for its support.}
\begin{abstract}
  In the first, mostly expository, part of this paper, a graded Lie
  algebra is associated to every group $G$ given with an $N$-series of
  subgroups. The asymptotics of the Poincar\'e series of this algebra
  give estimates on the growth of the group $G$. This establishes the
  existence of a gap between polynomial growth and growth of type
  $e^{\sqrt n}$ in the class of residually--$p$ groups, and gives
  examples of finitely generated $p$--groups of uniformly exponential
  growth.
  
  In the second part, we produce two examples of groups of finite
  width and describe their Lie algebras, introducing a notion of
  \emph{Cayley graph} for graded Lie algebras. We compute explicitly
  their lower central and dimensional series, and outline a general
  method applicable to some other groups from the class of branch
  groups.

  These examples produce counterexamples to a conjecture on the
  structure of just-infinite groups of finite width.
\end{abstract}
\maketitle
\thispagestyle{empty}

\section{Introduction}
The main goal of this paper is to present new examples of groups of
finite width and to give a method of proving that some groups from
the class of branch groups have finite width. This provides examples
of groups of finite width with a completely new origin and answers a
question asked by several mathematicians. We also give new examples of
Lie algebras of finite width associated to the groups mentioned above.

The first group we study, $\gr$, was constructed
in~\cite{grigorchuk:burnside} where it was shown to be an
infinite torsion group; later in~\cite{grigorchuk:gdegree} it was
shown to be of intermediate growth. The second group, $\grt$, was
already considered by the second author in 1979, but was rejected at
that time for not being periodic. We now know that it also has
intermediate growth~\cite{bartholdi-g:parabolic} and finite width.

Our interest in the finite width property comes from the theory of
growth of groups. Another important area connected to this property is
the theory of finite $p$-groups and the theory of pro-$p$-groups;
see~\cite{shalev:pro-p}, \cite[\S8]{shalev:finite-p}
and~\cite{klass-lg-p:fw} with its bibliography. More precisely, the
following was discussed by many mathematicians and stated by Zel'manov in Castelvecchio in 1996~\cite{zelmanov:castelvecchio}:
\begin{conj}\label{conj:zelmanov}
  Let $G$ be a just-infinite pro-$p$-group of finite width. Then
  $G$ is either solvable, $p$-adic analytic, or commensurable to a
  positive part of a loop group or to the Nottingham group.
\end{conj}
Our computations disprove this conjecture by providing a
counter-example, the profinite completion of $\gr$ (it is a
pro-$p$-group with $p=2$). Note that it exhibits a behaviour specific
to positive characteristic: indeed it was proved by Martinez and
Zel'manov in~\cite{zelmanov-m:nil} that unipotence and finite width
imply local nilpotence.

Before we give the definition of a group of finite width, let us
recall a classical construction of Magnus~\cite{magnus:lie},
described for instance in~\cite[Chapter~VIII]{huppert-b:fg2}. Given a
group $G$ and $\{G_n\}_{n=1}^\infty$ an $N$-series (i.e.\ a series of
normal subgroups with $G_1=G$, $G_{n+1}\le G_n$ and $[G_m,G_n]\leq
G_{m+n}$ for all $m,n\ge1$), there is a canonical way of associating
to $G$ a graded Lie ring
\begin{equation}\label{eq:liering}
  \Lie(G) = \bigoplus_{n=1}^\infty L_n,
\end{equation}
where $L_n=G_n/G_{n+1}$ and the bracket operation is induced by
commutation in $G$. Possible examples of $N$-series are the lower
central series $\{\gamma_n(G)\}_{n=1}^\infty$; for an integer $p$, the
\emph{lower $p$-central series} given by $P_1(G)=G$ and
$P_{n+1}(G)=P_n(G)^p[P_n(G),G]$; and, for a field $\Bbbk$, the series
of \emph{$\Bbbk$-dimension subgroups} $\{G_n\}_{n=1}^\infty$ defined
by
$$G_n = \{g\in G|\,g-1\in\Delta^n\},\qquad n=1,2,\dots$$
where $\Delta$ is the augmentation (or fundamental) ideal of the group
algebra $\Bbbk[G]$.

Tensoring the $\Z$-modules $L_n$ with a suitable field $\Bbbk$, we
obtain in~(\ref{eq:liering}) a graded Lie algebra $\Lie_\Bbbk(G)$. In
case the $N$-series chosen satisfies the additional condition
$G_n^p\le G_{pn}$, and $\Bbbk$ is a field of characteristic $p$, the
algebra $\Lie_\Bbbk(G)$ will then be a $p$-algebra or \emph{restricted
  algebra}; see~\cite{jacobson:restr} or~\cite[Chapter~V]{jacobson:lie},
the Frobenius operation on $\Lie_\Bbbk(G)$ being induced by raising to
the power $p$ in $G$. In this case the quotients $G_n/G_{n+1}$ are
elementary $p$-groups.

Many properties of a group are reflected in properties of its
corresponding Lie algebra. For instance, one of the most important
results obtained using the Lie method is the theorem of
Zel'manov~\cite{zelmanov:liemethod} asserting that if the Lie algebra
$\Lie_\Fp(G)$ associated to the dimension subgroups of a finitely
generated periodic residually-$p$ group $G$ satisfies a polynomial
identity then the group $G$ is finite ($\Fp$ is the prime field of
characteristic $p$). This result gives in fact a positive solution to
the Restricted Burnside
Problem~\cite{vaughan-lee-z:bp1,zelmanov:morebp,vaughan-lee-z:bp2,zelmanov:bp}.
Another example is the criterion of analyticity of pro-$p$-groups
discovered by Lazard~\cite{lazard:padanal}.

The Lie method also applies to the theory of growth of groups, as was
first observed in~\cite{grigorchuk:hp}. There the second author
proved that in the class of residually-$p$ groups there is a gap
between polynomial growth and growth of type $e^{\sqrt n}$. This
result was then generalized in~\cite[Theorem~D]{lubotzky-m:subgp} to the class
of residually-nilpotent groups, and in~\cite{ceccherini-g:unifexpo}
the Lie method was also used to prove that certain one-relator groups
with exponential-growth Lie algebra $\Lie_\Bbbk(G)$ have uniformly
exponential growth. If a group $G$ is finitely generated, then its Lie
algebra $\Lie_\Bbbk(G)=\bigoplus L_n\otimes\Bbbk$ is also finitely
generated, and the growth of $\Lie_\Bbbk(G)$ is by definition the
growth of the sequence $\{b_n = \dim(L_n\otimes\Bbbk)\}_{n=1}^\infty$.

The investigation of the growth of graded algebras related to groups
has its own interest and is related to other topics.  One of the first
results in this direction is the Golod-Shafarevich
inequality~\cite{golod-s:orig} which plays an
important role in group, number and field theories. The idea of Golod
and Shafarevich was used by Lazard in the proof of the aforementioned
criterion of analyticity (he even used the notation `gosha' for the
growth of the algebras).  Vershik and Kaimanovich
observed the relation between the growth of gosha, amenability, and
asymptotic behaviour of random walks (see Section~\ref{sec:amen} below).

For our purposes it will be sufficient to consider only the fields
$\Q$ and $\Fp$. Let $G_n$ be the corresponding series of dimension
subgroups, which is also an $N$-series, and let $\Lie_\Bbbk(G)$ be the
associated Lie algebra. If $\Lie_\Bbbk(G)$ is of polynomial growth of
degree $d\ge 0$, then the growth of $G$ is at least
$e^{n^{1-1/(d+2)}}$, and if $\Lie_\Bbbk(G)$ is of exponential growth,
then $G$ is of uniformly exponential growth.

If $\Bbbk=\Q$ and $G$ is residually-nilpotent and $b_n=0$ for some
$n$, then $G$ is nilpotent; indeed $G_n$ must be finite for that $n$,
whence $\gamma_n(G)$ is finite too, and since
$\bigcap_{k\ge1}\gamma_k(G)=1$ this implies that $\gamma_N(G)=1$ for
some $N$. It follows that $G$ has polynomial
growth~\cite{milnor:solvable}. In fact polynomial growth is equivalent
to virtual nilpotence~\cite{gromov:nilpotent}.

If $\Bbbk=\Fp$ and $G$ is a residually-$p$ group and $b_n=0$ for some
$n$, then $G$ is a linear group over a field, by Lazard's
theorem~\cite{lazard:padanal} and therefore has either polynomial
or exponential growth, by the Tits alternative~\cite{tits:linear}.

Finally, if $b_n\ge1$ for all $n$ then, independent of $\Bbbk$, the
growth of $G$ is at least $e^{\sqrt n}$. Keeping in mind that
polynomial growth $b_n\sim n^d$ of $\Lie_\Bbbk(G)$ implies a lower
bound $e^{n^{1-1/(d+2)}}$ for the growth of $G$, we conclude that
examples of groups with growth exactly $e^{\sqrt n}$ are to be found
amongst the class of groups for which the sequence
$\{b_n\}_{n=1}^\infty$ is uniformly bounded, or at least bounded in
average. This key observation leads to the notion of \emph{groups of
  finite width}. We present two versions of the definition:
\begin{defn}\label{defn:fw}
  Let $G$ be a group and $\Bbbk\in\{\Q,\Fp\}$ a field. If $\Bbbk=\Q$,
  assume $G$ is residually-nilpotent; if $\Bbbk=\Fp$, assume $G$ is
  residually-$p$.
  \begin{enumerate}
  \item $G$ has \emdef{finite $C$-width} if there is a constant $K$
    with $[\gamma_n(G):\gamma_{n+1}(G)]\le K$ for all $n$.
  \item $G$ has \emdef{finite $D$-width with respect to $\Bbbk$} if
    there is a constant $K$ with $b_n\le K$ for all $n$, where
    $\{b_n\}_{n=1}^\infty$ is the growth of $\Lie_\Bbbk(G)$
    constructed from the dimension subgroups.
  \end{enumerate}
\end{defn}
A third notion can be defined, that of \emdef{finite averaged width};
see~\cite{grigorchuk:hp}
or~\cite[Definition~I.1.ii]{klass-lg-p:fw}.  From our point of view
$D$-width is more natural; but the first notion is more commonly used
in the theory of finite $p$-groups and pro-$p$-groups, see for
instance~\cite[Definition~I.1.i]{klass-lg-p:fw}.
The examples we will produce are of finite width according to both
definitions. That one of our groups has finite width was conjectured
in~\cite{grigorchuk:hp}; it was proven that the numbers $b_n$ are
bounded in average. Rozhkov then confirmed this conjecture
in~\cite{rozhkov:lcs} by computing explicitly the $b_n$; but the proof
had gaps, one of which was filled in~\cite{rozhkov:habilitation}. We
fix another gap in the ``Technical Lemma~4.3.2''
of~\cite{rozhkov:habilitation} while simplifying and clarifying
Rozhkov's proof, and also outline a general method, connected to ideas
of Kaloujnine~\cite{kaloujnine:lcs}.

We recall in the next section known notions on algebras associated to
groups, and construct in Section~\ref{sec:torsue} a torsion group of
uniformly exponential growth. Section~\ref{sec:gptree} describes a
class of groups acting on rooted trees, and the next two sections
detail for two specific examples the indices of the lower central and
dimensional series. More specifically, we compute in
Theorem~\ref{thm:indG} and~\ref{thm:indS} the indices of these series
for the group $\gr$ and an overgroup $\grt$. We also obtain
in the process the structure of the Lie algebras $L(\gr)$ (associated
to the lower central series) and $\Lie_{\F_2}(\gr)$ (associated to the
dimension series) in Theorem~\ref{thm:lieG}, and that of $L(\grt)$
and $\Lie_{\F_2}(\grt)$ in Theorem~\ref{thm:lieS}. They are described
using Cayley graphs of Lie algebras, see
Subsection~\ref{subs:cayG}.

Throughout this paper groups shall act on the left. We use the
notational conventions $[x,y]=xyx^{-1}y^{-1}$ and $x^y=yxy^{-1}$.

Both authors wish to thank Aner Shalev and Efim Zelmanov for their
interest and generous contribution through discussions.

\section{Growth of Groups and Associated Graded Algebras}
Let $G$ be a group, $\{\gamma_n(G)\}_{n=1}^\infty$ the lower central
series of $G$, $\Bbbk\in\{\Q,\Fp\}$ a prime field,
$\Delta=\ker(\varepsilon)<\Bbbk[G]$ the augmentation ideal, where
$\varepsilon(\sum k_i g_i)=\sum k_i$ is the augmentation map
$\Bbbk[G]\to\Bbbk$, and $\{G_n\}_{n=1}^\infty$ the series of dimension
subgroups of $G$~\cite{zassenhaus:ordnen,jennings:gpring}. Recall that
$$G_n = \{g\in G|\,g-1\in\Delta^n\}.$$
The restrictions we impose on $\Bbbk$ are not important, as $G_n$
depends only on the characteristic of $\Bbbk$.  We suppose throughout
that $G$ is residually-nilpotent if $\Bbbk=\Q$ and is residually-$p$
if $\Bbbk=\Fp$.

If $\Bbbk=\Q$, then $G_n$ is the isolator of $\gamma_n(G)$, as was
proved in~\cite{jennings:gpringnilp} (see
also~\cite[Theorem~11.1.10]{passman:gr}
or~\cite[Theorem~IV.1.5]{passi:gr}); i.e.\ 
$$G_n = \sqrt{\gamma_n(G)} = \{g\in G|\,g^\ell\in\gamma_n(G)\text{ for
  an }\ell\in\N\}.$$
Note that in~\cite{passman:gr} these results are stated for finite
$p$-groups. They nevertheless hold in the more general setting of
residually-nilpotent or residually-$p$ groups.

If $\Bbbk=\Fp$, then $\gamma_n(G)\le G_n\le\sqrt{\gamma_n(G)}$, and
the $G_n$ can be defined in several different ways, for instance by
the relation
$$G_n = \prod_{i\cdot p^j\ge n}\gamma_i^{p^j}(G)$$
due to Lazard~\cite{lazard:nilp}, or recursively as
\begin{equation}\label{eq:recj}
  G_n = [G,G_{n-1}]G_{\lceil n/p \rceil}^p,
\end{equation}
where $\lceil n/p \rceil$ is the least integer greater than or equal
to $n/p$. In characteristic $p$, the series $\{G_n\}_{n=1}^\infty$ is
called the lower $p$-central, Brauer, Jennings, Lazard or Zassenhaus
series of $G$. The quotients $G_n/G_{n+1}$ are elementary abelian
$p$-groups and define the fastest-decreasing central series with the
property $G_n^p\le G_{np}$~\cite{jennings:gpringnilp}.

Let $$\mathcal A(G) = \mathcal A_\Bbbk(G) = \bigoplus_{n=0}^\infty
\Delta^n/\Delta^{n+1}$$ be the associative graded algebra with product
induced linearly from the group product
(see~\cite{passman:gr,passi:gr} for more details).

If $\Bbbk=\Q$, consider the following graded Lie algebras over $\Bbbk$:
$$\Lie(G)=\bigoplus_{n=1}^\infty\big(G_n/G_{n+1}\otimes_\Z\Q\big),\qquad
L(G)=\bigoplus_{n=1}^\infty\big(\gamma_n(G)/\gamma_{n+1}(G)\otimes_\Z\Q\big).$$
If $\Bbbk=\Fp$, consider the restricted Lie $\Fp$-algebra
$$\Lie_p(G)=\bigoplus_{n=1}^\infty\big(G_n/G_{n+1}\big).$$

Then Quillen's Theorem~\cite{quillen:ab} asserts that $\mathcal A(G)$ is
the universal enveloping algebra of $\Lie(G)$ in characteristic $0$
and is the universal $p$-enveloping algebra of $\Lie_p(G)$ in positive
characteristic.

Let us introduce the following numbers:
$$a_n(G) = \dim_\Bbbk(\Delta^n/\Delta^{n+1}),\qquad b_n(G) =
\rank(G_n/G_{n+1}).$$
Here by the rank of the $G$-module $M$ we mean
the torsion-free rank\\
 $\dim_\Q(M~\otimes~\Q)$ in characteristic $0$ and
the $p$-group rank $\dim_\Fp(M~\otimes~\Fp)$, equal to the minimal
number of generators, in positive characteristic. Note that in
zero-characteristic $b_n = \rank(\gamma_n(G)/\gamma_{n+1}(G))$,
because the natural map 
$$\gamma_n(G)/\gamma_{n+1}(G)\to G_n/G_{n+1}$$
has finite kernel and cokernel.

The following result is due to Jennings. The case $\Bbbk=\Fp$ appears
in~\cite{jennings:gpring} and the case $\Bbbk=\Q$ appears
in~\cite{jennings:gpringnilp}; but see also~\cite[Theorem~3.3.6
and~3.4.10]{passman:gr}.
\begin{equation}\label{eq:quillen}
  \sum_{n=0}^\infty a_n(G)t^n=\begin{cases}
    \prod_{n=1}^\infty\left(\frac{1-t^{pn}}{1-t^n}\right)^{b_n(G)}&\text{ if }\Bbbk=\Fp,\\
    \prod_{n=1}^\infty\left(\frac1{1-t^n}\right)^{b_n(G)}&\text{ if }\Bbbk=\Q.
  \end{cases}
\end{equation}

The series $\sum_{n=0}^\infty a_n(G)t^n$ is the Hilbert-Poincar\'e
series of the graded algebra $\mathcal A(G)$. The
equation~(\ref{eq:quillen}) expresses this series in terms of the
numbers $b_n(G)$; the relation between the sequences
$\{a_n(G)\}_{n=0}^\infty$ and $\{b_n(G)\}_{n=1}^\infty$ is quite
complicated. We shall be interested in asymptotic growth of series,
in the following sense:
\begin{defn}\label{defn:equiv}
  Let $f$ and $g$ be two functions $\R_+\to\R_+$. We write $f\precsim
  g$ if there is a constant $C>0$ such that $f(x)\le C+Cg(Cx+C)$ for
  all $x\in\R_+$, and write $f\sim g$ if $f\precsim g$ and $g\precsim
  f$.

  A series $\{a_n\}_{n=0}^\infty$ defines a function $f:\R_+\to\R_+$
  by $f(x)=a_{\lfloor x\rfloor}$, and for two series $a=\{a_n\}$ and
  $b=\{b_n\}$ we write $a\precsim b$ and $a\sim b$ when the same
  relations hold for their associated functions.
\end{defn}
The main facts are presented in the following statement:
\begin{prop}\label{prop:growth}
  Let $\{a_n\}$ and $\{b_n\}$ be connected by the
  one of the relations~(\ref{eq:quillen}). Then
  \begin{enumerate}
  \item $\{b_n\}$ grows exponentially if and only if $\{a_n\}$ does, and
    we have
    $$\limsup_{n\to\infty}\frac{\ln a_n}n = \limsup_{n\to\infty}\frac{\ln b_n}n.$$\label{itemg:2}
  \item If $b_n\sim n^d$ then $a_n\sim e^{n^{(d+1)/(d+2)}}$.\label{itemg:1}
  \end{enumerate}
\end{prop}
\begin{proof}
  We first suppose $\Bbbk=\Q$, and prove Part~\ref{itemg:2}
  following~\cite{berezny:subexp}. Let
  $A=\limsup( \ln a_n)/n$ and $B=\limsup (\ln b_n)/n$. Clearly
  $A\ge B$ as $a_n\ge b_n$ for all $n$; we now prove that $A\le B$. Define
  $$f(z)=\prod_{n=1}^\infty(1-e^{-nz})^{-b_n},$$
  viewed as a complex analytic function in the half-plane $\Re(z)>B$.
  We have $|1-e^{-nz}|^{-1}\le(1-e^{-n\Re z})^{-1}$, from which $|f(z)|\le
  f(\Re z)$. Now applying the Cauchy residue formula,
  $$a_n=\frac1{2\pi}\int_{-\pi}^\pi f(u+iv)e^{n(u+iv)}dv\le\frac1{2\pi}\int_{-\pi}^\pi |f(u+iv)|e^{nu}dv\le e^{nu}f(u)$$
  for all $u>B$, so
  $$A=\limsup_{n\to\infty}\frac{\ln a_n}n\le\limsup_{u>B,n\to\infty}\left(u+\frac{\ln f(u)}n\right)=B.$$
  For $\Bbbk=\Fp$, Part~\ref{itemg:2} holds \emph{a fortiori}.
  
  Part~\ref{itemg:1} for $\Bbbk=\Q$ is a consequence of a result by
  Meinardus~(\cite{meinardus:asympt}; see
  also~\cite[Theorem~6.2]{andrews:part}).  More precisely, when
  $b_n=n^d$, his result implies that
  $$a_n\approx\frac{e^{\zeta'(-d)}}{\sqrt{2\pi(d+2)n}}\left(\frac{(d+1)!\zeta(d+2)}{n}\right)^{\frac{1-2\zeta(-d)}{2+4d}}e^{n\frac{d+2}{d+1}\left(\frac{(d+1)!\zeta(d+2)}{n}\right)^{\frac1{1+2d}}},$$
  where `$\approx$' means that the quotient tends to $1$ as
  $n\to\infty$, and $\zeta$ is the Riemann zeta function.
  
  We sketch the proof for $\Bbbk=\Q$ below: we suppose that $b_n\sim
  n^d$, so $A=B=0$ by Part~\ref{itemg:2}, and compute
\begin{eqnarray*}
 \frac{d}{du}\ln f(u)=\sum_{n=1}^\infty-b_n\frac{-ne^{-nu}}{1-e^{-nu}}
  & \sim & \frac1{u^{d+2}}\sum_{n=1}^\infty\frac{(nu)^{d+1}}{e^{nu}-1}u\\
  &\sim& \frac1{u^{d+2}}\int_0^\infty\frac{w^{d+1}}{e^w-1}dw=\frac C{u^{d+2}}.
\end{eqnarray*}
  Thus $\ln f(u)\sim C/u^{d+1}$, and the inequality
  $$\log a_n\le nu+\log f(u)\sim nu+C/u^{d+1}$$
  is tight by the saddle-point principle when the right-hand side is
  minimized. This is done by choosing $u=n^{-1/(d+2)}$, whence as
  claimed $\log a_n\sim n^{1-1/(d+2)}$.

  Finally, we show that~(\ref{eq:quillen}) yields the same asymptotics
  when $\Bbbk=\Fp$ as when $\Bbbk=\Q$. Clearly
  $$\prod_{n=1}^\infty(1+t^n)^{b_n}\le\prod_{n=1}^\infty(1+t^n+\dots+t^{(p-1)n})^{b_n}\le\prod_{n=1}^\infty(1+t^n+\dots)^{b_n}$$
  for all $p\ge2$, where for two power series $\sum e_t^n$ and $\sum
  f_nt^n$ the inequality $\sum e_t^n\le\sum f_nt^n$ means that $e_n\le
  f_n$ for all $n$. It thus suffices to consider the case $p=2$. For
  this purpose define
  $$g(z)=\prod_{n=1}^\infty(1+e^{-nz})^{b_n},$$
  and compare the series developments of $\log(f)$ and $\log(g)$ in
  $e^{-z}$. From $-\log(1-z)=\sum_{n\ge1}\frac{z^n}n$ it follows that
  \begin{gather*}
    \log f(z)=\sum_{n\ge1}f_ne^{-nz},\quad f_n=\sum_{d|n}\frac1d,\\
    \log g(z)=\sum_{n\ge1}g_ne^{-nz},\quad g_n=\sum_{d|n}\frac{(-1)^{d+1}}d,
  \end{gather*}
  so both series have the same odd-degree coefficients, and thus $\log
  f\sim\log g$. Their exponentials then have the same asymptotics; more
  precisely, $f_n\le g_{2n-1}$ for all $n$, so $e^z\log f(2z)\le\log
  g(z)$ termwise, and $f(2z)\le g(z)$.
\end{proof}


\subsection{Growth of Groups} Let $G$ be a finitely
generated group with a fixed semigroup system $S$ of generators (i.e.\ 
such that every element $g\in G$ can be expressed a product
$g=s_1\dots s_n$ for some $s_i\in S$). Let $\gamma_G^S(n)$ be the
growth function of $(G,S)$; recall that it is
$$\gamma_G^S(n) = \#\{g\in G|\,|g|\le n\},$$
where $|g|$ is the minimal number of generators required to express
$g$ as a product.

The following observations are well-known:
\begin{lem}\label{lem:equiv}
  Let $G$ be a group and consider two finite generating sets $S$ and
  $T$. Then $\gamma_G^S\sim\gamma_G^T$, with $\sim$ given in
  Definition~\ref{defn:equiv}.
\end{lem}
It is then meaningful to consider the \emdef{growth} $\gamma_G$ of
$G$, which is the $\sim$-equivalence class containing its growth
functions $\gamma_G^S$.

\begin{lem}\label{lem:subgrowth}
  Let $G$ be a finitely generated group, $H<G$ a finitely
  generated subgroup and $K$ a quotient of $G$. Then 
  $\gamma_H\precsim\gamma_G$ and $\gamma_K\precsim\gamma_G$.
\end{lem}
\begin{proof}
  Let $S$ be a finite generating set for $H$; choose a generating set
  $T\supset S$ for $G$. Apply Definition~\ref{defn:equiv} with $C=1$
  to obtain $\gamma_H^S\precsim\gamma_G^T$. Clearly
  $\gamma_K^T(n)\le\gamma_G^T(n)$ for all $n$.
\end{proof}

\begin{lem}[\cite{grigorchuk:hp}]\label{lem:alggp}
  For any field $\Bbbk$ and any group $G$ with generating set $S$ the
  inequalities $a_n(G)\le\gamma_G^S(n)$ hold for all $n\ge0$.
\end{lem}
\begin{proof}
  Fix a generating set $S$. The identities
  $$xy-1=(x-1)+(y-1)+(x-1)(y-1),\qquad x^{-1}-1=-(x-1)-(x-1)(x^{-1}-1)$$
  show that
  $$xy-1\equiv(x-1)+(y-1),\qquad x^{-1}-1\equiv-(x-1)\mod\Delta^2,$$
  so $\Delta^n$ is generated over $\Bbbk$ by $\Delta^{n+1}$ and
  elements of the form
  $$x_0(s_1-1)x_1(s_2-1)\dots(s_n-1)x_n,$$
  for all $s_i\in S$ and $x_i\in\Bbbk[G]$. Now
  $x_i\equiv\varepsilon(x_i)\in\Bbbk$ modulo $\Delta$, so
  $\Delta^n/\Delta^{n+1}$ is spanned by the
  $$(s_1-1)(s_2-1)\dots(s_n-1),\qquad s_i\in S.$$
  All these elements are in the vector subspace $S_n$ of $\Bbbk[G]$
  spanned by products of at most $n$ generators, and by definition
  $S_n$ is of dimension $\gamma_G^S(n)$.
\end{proof}

\begin{cor}\label{cor:alggp}
  $\{a_n(G)\}_{n=0}^\infty\precsim\gamma_G$.
\end{cor}

Combining Proposition~\ref{prop:growth} and Lemma~\ref{lem:alggp}, we
obtain as
\begin{cor}\label{cor:superrad}
  If there exist $C>0$ and $d\ge0$ such that $b_n\ge Cn^d$ for all
  $n$, then $\gamma_G(n)\succsim e^{1-1/(d+2)}$. In particular, if
  $b_n\neq0$ for all $n$, then $\gamma_G(n)\succsim e^{\sqrt n}$.
\end{cor}

We shall say a group $G$ is of \emdef{subradical growth} if
$\gamma_G\precnsim e^{\sqrt n}$.
\begin{thm}[\cite{grigorchuk:hp}]\label{thm:alagromov}
  Let $G$ be a finitely generated residually-$p$ group. If $G$ is of
  subradical growth then $G$ is virtually nilpotent and
  $\gamma_G(n)\sim n^d$ for some $d\in\N$.
\end{thm}
\begin{proof}
  By the previous corollary, $b_n(G)=0$ for some $n$. Consider the
  $p$-completion $\widehat G$ of $G$. As Lie algebras, $\Lie_\Fp(G)$
  and $\Lie_\Fp(\widehat G)$ coincide, so $b_n(\widehat G)=0$. By
  Lazard's criterion $\widehat G$ is an analytic
  pro-$p$-group~\cite{lazard:padanal} and thus is linear over a field.
  Since $G$ is residually-$p$ it embeds in $\widehat G$ so is also
  linear. By the Tits alternative~\cite{tits:linear} either $G$ contains
  a free group on two generators (contradicting the assumption on the
  growth of $G$) or $G$ is virtually solvable. By the results of Milnor and  Wolf every virtually solvable group is either of
  exponential growth or is virtually
  nilpotent~\cite{milnor:solvable,wolf:solvable}. The asymptotic
  growth is invariant under taking finite-index subgroups, and the
  growth of a nilpotent group is polynomial of degree
  $\sum_{k\ge1}kb_k$, as was shown by Guivarc'h and  Bass~\cite{guivarch:poly1,guivarch:poly2,bass:nilpotent}.
\end{proof}

In the class of residually-$p$ groups, Theorem~\ref{thm:alagromov}
improves Gromov's result~\cite{gromov:nilpotent} that a
finitely generated group $G$ having polynomial growth is virtually
nilpotent, in that the assumption is weakened from `polynomial growth'
to `subradical growth'. Lubotzky and Mann have shown the
same result for residually nilpotent groups of subradical growth. It
is not known whether subradical growth does imply virtual nilpotence,
and whether there exist groups of precisely radical growth. Certainly
the right place to look for such examples is among groups of finite
width, or groups satisfying some tight condition on the growth of
their $b_n$.

Therefore new examples of groups of finite width are of special
interest. Below we shall give two examples of such groups and outline
a method of constructing new examples; but first a consequence
of~\ref{thm:alagromov} is
\begin{thm}\label{thm:superrad}
  The growth $\gamma_\gr$ of the group $\gr$ satisfies
  $$e^{\sqrt n}\precsim\gamma_\gr(n)\precsim e^{n^{1/(1-\log_2\eta)}},$$
  where $\eta$ is the real root of $X^3+X^2+X-2$.
\end{thm}
\begin{proof}  
  If $\gr$ were nilpotent it would be finite, as it is finitely
  generated and torsion; since it is infinite~\ref{thm:alagromov}
  yields the left inequality.

  The right inequality was proven by the first author
  in~\cite{bartholdi:upperbd}, using purely combinatorial techniques.
\end{proof}
Note that the estimate from below can be obtained directly as
in~\cite{grigorchuk:gdegree}, by showing that for an appropriate
$S$ the growth function $\gamma_G^S$ satisfies
$$\gamma_G^S(4n)\ge\gamma_G^S(n)^2.$$
The second author conjectured in 1984 that the left inequality is in
fact an equality, but Leonov recently announced that
this is not the case~\cite{leonov:lowerbd}.

For our second example $\grt$ it is only known that
$$e^{\sqrt n}\precsim\gamma_\grt\precnsim e^n,$$
as is shown in~\cite{bartholdi-g:parabolic}.

Lemma~\ref{lem:alggp} can also be used to study uniformly exponential
growth, as was observed in~\cite{ceccherini-g:unifexpo}. Let
$$\omega_G^S=\lim_{n\to\infty}\sqrt[n]{\gamma_G^S(n)}$$
be the base of exponential growth of $G$ with respect to the
generating set $S$ and let $\omega_G=\inf_S\omega_G^S$, the infimum
being taken over all finite generating sets.
\begin{defn}
  The group $G$ has \emdef{uniformly exponential growth} if
  $\omega_G>1$.
\end{defn}

\noindent
(See~\cite{gromov:metriques} for the original definition and
motivations, and~\cite{grigorchuk-h:groups} for more details on this
notion.) For instance, the free groups of rank $\ge2$, and more
generally, the non-elementary hyperbolic groups have uniformly
exponential growth~\cite{koubi:unifexpo}. It is currently not known
whether there exists a group of exponential but not uniformly
exponential growth.
\begin{cor}
  If for some $\Bbbk\in\{\Q,\Fp\}$ the algebra $\mathcal A_\Bbbk(G)$
  has exponential growth then $G$ has uniformly exponential growth.
  (We do not need here the assumption that $G$ is residually-$p$ or
  residually nilpotent.)
\end{cor}

In the next section we will combine this idea with the Golod-Shafarevich
construction to produce examples of finitely generated
residually finite $p$-groups of uniformly exponential growth.

\section{Torsion Groups of Uniformly Exponential Growth}\label{sec:torsue}
As a reference to the Golod-Shafarevich construction we recommend the
original paper~\cite{golod-s:orig}, one of the
books~\cite{herstein:rings,koch:galois},
or~\cite[\S~VIII.12]{huppert-b:fg2}.

Consider the free associative algebra $A$ over the field $\Fp$ on the
generators $x_1,\dots,x_d$ for some $d\ge2$. The algebra $A$ is
graded: $A=\bigoplus_{n=0}^\infty A_n$ where $A_n$ is spanned by the
monomials of degree $n$, with $A_0=\Fp1$. Elements of the subspace
$A_n$ are called \emph{homogeneous of degree $n$}.

Consider an ideal $\mathcal I$ in $A$ generated by $r_1$ homogeneous
elements of degree $1$, $r_2$ of degree $2$, etc. (We make this
homogeneity assumption for simplicity; it is not necessary, as was
indicated in~\cite{koch:galois}.) Let $B=A/\mathcal
I$. Then $B$ is also a graded algebra: $B=\bigoplus_{n=0}^\infty B_n$
and if $H_B(t)=\sum_{n=0}^\infty d_n t^n$ be the Hilbert-Poincar\'e
series of $B$, i.e.\ $d_n=\dim_\Fp B_n$, then the Golod-Shafarevich
inequality
\begin{equation}\label{eq:golod-s}
  H_B(t)(1-dt+H_R(t))\ge1
\end{equation}
holds; here $H_R(t)=\sum_{n=1}^\infty r_nt^n$, and for the comparison
of two power series the same agreement holds as in the previous section.

Suppose that for some $\xi\in(0,1)$ the series $H_R(t)$ converges
at $\xi$ and $1-d\xi+H_R(\xi)\le0$. Then the series $H_B(t)$ cannot
converge at $t=\xi$, so the coefficients $d_n$ of $H_B(t)$ grow
exponentially and
$$\limsup_{n\to\infty}\sqrt[n]{d_n}\ge\frac1\xi.$$

Golod proves in~\cite{golod:nil} that $\mathcal I$ can be chosen in such a
way that the ideal $\mathcal D=\bigoplus_{n=1}^\infty B_n$ will be a
$p$-nilalgebra (i.e.\ for all $y\in\mathcal D$ there is an $n\in\N$
such that $y^{p^n}=0$).

The construction of the relators goes as follows: enumerate first as
$\{y_k\}_{k=1}^\infty$ all elements of the algebra $A$ (this is
possible since $A$ is countable). Start with $\mathcal I_0=0$; then if
$y_k$ is not a nilelement of $A/\mathcal I_{k-1}$ take $\ell_k\ge3$
sufficiently large so that the least degree of monomials in
$y_k^{p^{\ell_k}}$ is larger than all degrees of monomials in
$\mathcal I_{k-1}$. Construct $\mathcal I_k$ by adding to $\mathcal
I_{k-1}$ all homogeneous parts of the polynomial $y_k^{p^{\ell_k}}$.
Let finally $\mathcal I=\bigcup_{n=0}^\infty\mathcal I_n$.

The numbers $r_k$ will then all be $0$ or $1$ with $r_k=0$ for
$k<p^3$, so taking $\xi=3/4$ we have
$$1-d\xi+H_R(\xi)\le1-2\xi+\frac{\xi^{2^3}}{1-\xi}<0$$
and $B=A/\mathcal I$ is of exponential growth at least $(4/3)^n$. Let
$\overline x_1,\dots,\overline x_d$ be the images of $x_1,\dots,x_d$
in $B$, and let $G$ be the group generated by the elements
$s_i=1+\overline x_i$; they are invertible because the $\overline x_i$
are $p$-nilelements and $B$ is of characteristic $p$. The vector
subspace of $B$ spanned by $G$ is $B$ itself, so $B$ is a quotient of
the group algebra $\Fp[G]$.
\begin{thm}
  $G$ is a finitely generated residually finite $p$-group of uniformly
  exponential growth.
\end{thm}
\begin{proof}
  That $G$ is a $p$-group was observed by Golod and follows from the
  fact that $\mathcal D=\bigoplus_{n=1}^\infty B_n$ is a
  $p$-nilalgebra. Let $\pi$ be the natural map $\Fp[G]\to B$. Then
  $\mathcal D$ is generated by $\pi(\Delta)$ and more generally
  $\bigoplus_{n=N}^\infty B_n=\pi(\Delta^N)$, so by
  Lemmata~\ref{lem:alggp} and~\ref{lem:subgrowth} there is a $\xi<1$
  such that the estimate
  $$\frac1{\xi^n}\le\dim_\Fp B_n\le a_n(G)\le\gamma_G^T(n),\qquad n=1,2,\dots$$
  holds for any system $T$ of generators of $G$.
\end{proof}

\section{Growth of Algebras and Amenability}\label{sec:amen}
As was mentioned in the introduction, there is an interesting question
(due to Vershik) on the relation between the amenability of a group
and the growth of related algebras. Let us formulate our version of
this question:

\newpage
\begin{problem}

{\rm 1.} Let $G$ be amenable. Does $b_n(G)$ grow subexponentially for
    any field $\Bbbk$?
\\
{\rm 2.} Suppose $G$ is residually nilpotent (or residually-$p$) and
    $b_n(G)$ grows subexponentially for the field $\Q$ (or $\Fp$). Is
    then the group $G$ amenable?

\end{problem}

There is a chance that for at least one of these questions the answer
is affirmative.

For solvable groups (which are amenable) the associated algebras have
subexponential growth, as follows from computations by
Petrogradski\u\i~\cite{petrogradsky:growth,petrogradsky:growth2};
his results are based on computations for free polynilpotent algebras
by Bokut'~\cite{bokut:basis}. See also
Egorychev~\cite{egorychev:intsum} and
Berezni\u\i~\cite{berezny:subexp} for partial results.

On the other hand there is some similarity between the asymptotics of
random walks on solvable groups and the growth of
$b_n(G)$~\cite{kaimanovich:rw} which gives a hope that subexponential
growth of algebras implies (under the residuality hypothesis)
subexponential decay of the probability of returning to the origin for
symmetric random walks on a group. Then Kesten's
criterion~\cite{kesten:rwalks} can be invoked to imply the
amenability of $G$.

\section{Groups Acting on Rooted Trees}\label{sec:gptree}
We now consider examples of groups whose lower central series and
dimension series we can compute explicitly. Let $\Sigma$ be a finite
alphabet, and $\Sigma^*$ the set of finite sequences over $\Sigma$.
This set has a natural rooted tree structure: the vertices are finite
sequences, and the edges are all the $(\sigma,\sigma s)$ for
$\sigma\in\Sigma^*$ and $s\in\Sigma$; the root vertex is $\emptyset$,
the empty sequence. By $\aut(\Sigma^*)$ we mean the bijections of
$\Sigma^*$ that preserve the tree structure, i.e.\ preserve length and
prefixes. We write $\sigma\Sigma^*$ for the subtree of $\Sigma^*$
below vertex $\sigma$: it is isomorphic to $\Sigma^*$ but rooted at $\sigma$.

Let $G$ be a finitely generated subgroup of $\aut(\Sigma^*)$ acting
transitively on $\Sigma^n$ for all $n$ (such an action will be called
\emdef{spherically transitive}.) We denote by $\stab_G(\sigma)$ the
stabilizer of the vertex $\sigma$ in $G$, and by $\stab_G(n)$ the
stabilizer of all vertices of length $n$.  An arbitrary element
$g\in\stab_G(n)$ can be identified with a tuple
$(g_\sigma)_{|\sigma|=n}$ of tree automorphisms; we write this
monomorphism
$$\phi_n:\stab_G(n)\hookrightarrow\prod_{\sigma\in\Sigma^n}\aut(\sigma\Sigma^*).$$
We define the \emdef{vertex group} or \emdef{rigid stabilizer}
$\rist_G(\sigma)$ of the vertex $\sigma$ by
$$\rist_G(\sigma) = \{g\in
G|\,g\tau=\tau\;\forall\tau\in\Sigma^*\setminus\sigma\Sigma^*\},$$
and the \emdef{$n^{th}$ rigid stabilizer} as the group generated by
the length-$n$ vertex groups: $\rist_G(n)=\langle
\rist_G(\sigma):\,|\sigma|=n\rangle$.  Since $G$ acts transitively on
$\Sigma^n$ the vertex groups of vertices at level $n$ are all
conjugate. Therefore $\rist_G(n)$ is a direct product of $|\Sigma|^n$
copies of $\rist_G(\sigma)$ for a $\sigma$ of length $n$.

\begin{defn}\label{defn:branch}
  A finitely generated group $G$ is called a \emdef{branch group} if
  \begin{enumerate}
  \item $G$ acts faithfully on $\Sigma^*$ and transitively on
    $\Sigma^n$ for all $n\ge0$;
  \item $[G:\rist_G(n)]$ is finite for all $n\ge0$.
  \end{enumerate}
\end{defn}

\subsection{The Modules $V_n$}\label{subs:modules}
Let $G$ be a group acting on a regular rooted tree $\Sigma^*$, where
$\Sigma$ contains $p$ elements for some prime $p$; for ease of
notation suppose $\Sigma=\Fp$. Assume moreover that at each vertex $G$
acts as a power of the cyclic permutation $\epsilon=(0,1,\dots,p-1)$ of
$\Sigma$. Let $V_n=\Fp[G/\stab_G(0^n)]$; it is a vector space of
dimension $p^n$, as $G$ acts transitively on $\Sigma^n$, and has a
natural $G$-module structure coming from the action of $G$ on
$G/\stab_G(0^n)$.  Identify $G/\stab_G(0^n)$ with the set $\Sigma^n$ of
vertices at level $n$, and also with the set of monomials over
$\{X_1,\dots,X_n\}$ of degree $<p$ in each variable, by
$$\sigma=\sigma_1\dots\sigma_n\leftrightarrow X_1^{\sigma_1}\dots X_n^{\sigma_n}.$$
Under this identification, we can write
\begin{eqnarray*}
V_n &=& \Fp[X_1,\dots,X_n]/(X_1^p-1,\dots,X_n^p-1)\\
&=&\Fp[X_1]/(X_1^p-1)\otimes\dots\otimes\Fp[X_n]/(X_n^p-1).
\end{eqnarray*}
We write $g\sigma$ the action of $g\in G$ on $\sigma\in V_n$, and
denote by $[g,\sigma]=\sigma-g\sigma$ the ``Lie action'' of $G$ on
$V_n$. For $r\in\{0,\dots,p^n-1\}$ we write $r=r_n\dots r_1$ in base $p$,
and define
\begin{gather*}
  v_n^r = (1-X_1)^{r_1}\dots(1-X_n)^{r_n}\in V_n,\\
  V_n^r = \langle v_n^r,\dots,v_n^{p^n-1}\rangle.
\end{gather*}
We extend the last definition to $V_n^r=0$ when $r\ge p^n$. There is a
natural projective sequence
$$\dots\to V_n\to V_{n-1}\to\dots\to V_0=\Fp$$
of $G$-modules, and at each step $n$ a sequence of $V_n$-submodules
$$V_n^{p^n}=0\subset\dots\subset V_n^{p^n-p^{n-1}}\subset\dots\subset
V_n^1\subset V_n^0=V_n$$
each having codimension $1$ in the next.
Moreover $V_{n-1}^{p^{n-1}-i}$ and $V_n^{p^n-i}$ are naturally
isomorphic under multiplication by $(1-X_n)^{p-1}$; thus
$V_n^{p^n-p^{n-1}}$ is isomorphic to $V_{n-1}^0=V_{n-1}$ as a $G$-module.

\begin{lem}\label{lem:dec1}
  
{\rm 1.} The inclusion $[G,V_n^r]\subset V_n^{r+1}$ holds for all
    $n$ and all $r$.
\\
{\rm 2.} If $G$ contains for all $m\le n$ an element $g_m$ such that
    $$g_m(0^m)=0^{m-1}1,\qquad g_m(\sigma x)=\sigma'x\quad\forall\sigma\in\Sigma^{m-1}\setminus\{0^{m-1}\},x\in\Sigma$$
    (where in the second condition $\sigma'$ is an arbitrary function
    of $\sigma$), then $[G,V_n^r]=V_n^{r+1}$ for all $n$ and all
    $r$.

\end{lem}
A $G$-module $V$ having the property $\dim V^{(n)}/V^{(n+1)}=1$ for
all $n$, where the $V^{(n)}$ are defined inductively by $V^{(0)}=V$
and $V^{(n+1)}=[G,V^{(n)}]$ is called \emdef{uniserial}. This notion
was introduced by Leedham-Green~\cite{leedham-green:fpg}; see
also~\cite[page~111]{dixon-:prop}.

Note that every element of $G$ can be described by a colouring
$\{g_\sigma\}_{\sigma\in\Sigma^*}$ of the vertices of $\Sigma^*$ by
elements of the cyclic group $C_p=\langle\epsilon\rangle$. 
The condition in the lemma amounts to the existence, for all $m$, of an
element $g_m$ whose colouring is $\epsilon$ at the vertex $0^m$, and is
$1$ on all other vertices of the $m$-th level as well as on all vertices
$0^i$, for $i<m$.
 Note also that this implies that
the action is spherically transitive.
\begin{proof}
  We proceed by induction on $(n,r)$ in lexicographic order. For $n=0$
  the claim holds trivially; suppose thus $n\ge1$. In order to prove
  $[G,V_n^r]\subset V_n^{r+1}$, it suffices to check that for all
  $g\in G$ we have $[g,v_n^r]\in V_n^{r+1}$, as the $V_n^r$ form an
  ascending tower of subspaces. During the proof we will consider
  $V_{n-1}$ as a subspace of $V_n$; beware though that it is not a
  submodule. We shall write `$\ast$' for the action of $G$ on
  $V_{n-1}\subset V_n$, and `$\cdot$' for that of $G$ on $V_n$.
  
  Observe that if $v\in V_{n-1}$ then $g\cdot(vX_n^i)=(g\cdot
  v)X_n^i$. Thus $g\cdot v-g\ast v$ is always divisible by $1-X_n$
  because if $g\cdot v=\sum_{s=0}^{p-1}\Psi_sX_n^s$ for some
  $\Psi_s\in V_{n-1}$ then $g\ast v=\sum_{s=0}^{p-1}\Psi_s$ and
  \begin{equation}\label{eq:tlemma}
    g\cdot v-g\ast v=(1-X_n)\sum_{s=1}^{p-1}-\Psi_s(1+X_n+\dots+X_n^{s-1}).
  \end{equation}
 Write $r=r_n\dots r_1$ in base $p$. For some $\Phi$ and $\Psi_s$ in
  $V_{n-1}$, we may write
  $$v_n^r=\Phi(1-X_n)^{r_n},\qquad g\cdot v_n^r=\sum_{s=0}^{p-1}\Psi_sX_n^s(1-X_n)^{r_n}.$$
  Then by induction
  $$[g,v_n^r] = \underbrace{\underbrace{\left(\Phi-\sum_{s=0}^{p-1}\Psi_s\right)}_{\in V_{n-1}^{(r+1)\modop p^{n-1}}}(1-X_n)^{r_n}}_{\in V_n^{r+1}}+\underbrace{\sum_{s=0}^{p-1}\Psi_s(1-X_n^s)(1-X_n)^{r_n}}_{\in V_n^{(r_n+1)p^{n-1}}\subseteq V_n^{r+1}},$$
  as in the second summand $(1-X_n^s)(1-X_n)^{r_n}$ is divisible by
  $(1-X_n)^{r_n+1}$. This proves the first claim of the lemma.

  Next, we prove $[G,V_n^r]\supset V_n^{r+1}$ by showing that
  $v_n^{r+1}\in[G,V_n^r]$. As above, write $r=r_n\dots r_1$ in base
  $p$. If $(r_1,\dots,r_{n-1})\neq (p-1,\dots,p-1)$, we have
  $v_n^{r+1}=v_{n-1}^{r+1\modop p^{n-1}}(1-X_n)^{r_n}$, and by induction
  $v_{n-1}^{r+1\modop p^{n-1}}=\sum_s\alpha_s[g_s,v_{n-1}^{i_s}]$ for
  some $\alpha_s\in\Fp$, $g_s\in G$ and $i_s\ge r\modop p^{n-1}$.
  Then
  $$v_n^{r+1}=\underbrace{\sum_s\alpha_s\left[g_s,v_n^{i_s+r_np^{n-1}}\right]}_{\in[G,V_n^r]}+\underbrace{\sum_s\alpha_s\left(g_s\cdot v_n^{i_s+r_np^{n-1}}-(g_s\ast v_{n-1}^{i_s})(1-X_n)^{r_n}\right)}_{\in V_n^{(r_n+1)p^{n-1}}\subseteq V_n^{r+2}\subseteq[G,V_n^{r+1}]\subseteq[G,V_n^r]}$$
  where the last inclusions hold by~(\ref{eq:tlemma}) and induction.
  Finally, if $r=(r_n+1)p^{n-1}-1$, note that
  \begin{align*}
    v_n^r &= (1+X_1+\dots+X_1^{p-1})\cdots(1+X_{n-1}+\dots+X_{n-1}^{p-1})(1-X_n)^{r_n}\\
    &= (1-X_n)^{r_n} + P(1-X_n)^{r_n},
  \end{align*}
  where $P=\sum_{\sigma\in\Sigma^{n-1}\setminus\{0^{n-1}\}}X_1^{\sigma_1}\cdots
  X_n^{\sigma_n}$ is invariant under $g_n$; thus
  $$v_n^{r+1} = (1-X_n)^{r_n} - X_n(1-X_n)^{r_n} = [g_n,v_n^r]\in[G,V_n^r].$$
\end{proof}

The strategy we follow to compute the lower central series or
dimension series of $G$ in the examples of Sections~\ref{sec:G}
and~\ref{sec:S} is the following:
\begin{itemize}
\item We recognize some $\gamma_m(G)$ or $G_m$ as a subgroup of $G$
  simply obtained from rigid stabilizers in $G$.
\item We identify a quotient $\gamma_m(G)/N$ or $G_m/N$ with a direct
  sum of copies of the module $V_n$ defined above, for an appropriate
  subgroup $N$.
\item We show that $N$ is a further term of the lower central or
  dimensional series, allowing the process to repeat.
\end{itemize}
Then the exact terms of the lower central or dimension series are
obtained by pulling back the appropriate $V_n^r$ through the
identification.

\section{The Group $\gr$}\label{sec:G}
Let $\Sigma=\F_2$, the field on two elements. For $x\in\F_2$ set
$\overline x=1-x$, and define the automorphisms $a$, $b$, $c$, $d$ of
$\Sigma^*$ as follows:
$$a(x\sigma)=\overline x\sigma,$$
\begin{alignat*}{2}
  b(0x\sigma)&=0\overline x\sigma,&\qquad b(1\sigma)&=1c(\sigma),\\
  c(0x\sigma)&=0\overline x\sigma,&\qquad c(1\sigma)&=1d(\sigma),\\
  d(0x\sigma)&=0x\sigma,&\qquad d(1\sigma)&=1b(\sigma).
\end{alignat*}
Thus for instance $b$ acts on the subtree $0\Sigma^*$ as $c$, while
$c$ acts on it as $d$, etc. Note that all generators are of order $2$
and $\{1,b,c,d\}$ forms a Klein group. Set $\gr=\langle
a,b,c,d\rangle$. For ease of notation, we shall identify elements of
$\stab_\gr(n)$ with their image under $\phi_n$ by writing
$\phi_n(g)=(g_1,\dots,g_{2^n})_n$ (omitting the subscript $n$ if it is
obvious from context); for instance we will write $b=(a,c)$, $c=(a,d)$
and $d=(1,b)$.
Set $x=[a,b]$, and set
$$K=\langle x\rangle^\gr = \langle x,(x,1),(1,x)\rangle.$$
Note that
$(x,1)=[b,d^a]$ and $(1,x)=[b^a,d]$. Also, $K$ is a subgroup of finite
index (actually index $16$) in $\gr$, and contains $K\times K$ as a
subgroup of finite index (actually index $4$); for more details
see~\cite{harpe:cgt} or~\cite{bartholdi-g:parabolic}. Set also
$T=\langle x^2\rangle^\gr=K^2$, and for any $Q\le K$ define
$Q_m=Q\times\dots\times Q$ ($2^m$ copies). Clearly $Q_m\le\stab_\gr(m)$
and acts on each subtree starting on level $m$ by the corresponding
factor. For $m\ge1$ set $N_m=K_m\cdot T_{m-1}$.

For $m\ge2$, we have $\rist_\gr(m)=K_{m-2}$, so $\gr$ is a branch group.

\begin{lem}\label{lem:isoG} The mapping
  $$\alpha\oplus\beta:N_m/N_{m+1}\longrightarrow V_m\oplus V_{m-1}$$
  is an isomorphism for all $m$, where the $V_m$ are the modules defined
  in Subsection~\ref{subs:modules}, $\alpha$ maps
  $(1,\dots,1,x,1,\dots,1)\in K_m$ to the monomial in $V_m$
  corresponding to the vertex at the $x$'s position,
  and $\beta$ maps $(1,\dots,1,x^2,1,\dots,1)\in T_{m-1}$ to the
  corresponding monomial in $V_{m-1}$.
\end{lem}
\begin{proof}
  We first suppose $m=1$. Then $N_1/N_2=\langle
  x^2,(1,x),(x,1)\rangle/N_2$; it is easy to check that
  $x^4=(x^2,x^2)$ modulo $K_2$, so all generators of $N_1/N_2$ are of
  order $2$. Further, $[x^2,(1,x)]\in K_2$ and $[x^2,(x,1)]\in K_2$,
  so the quotient $N_1/N_2$ is the elementary abelian group $2^3$, and
  $\alpha\oplus\beta$ is an isomorphism in that case.

  For $m>1$ it suffices to note that both sides of the isomorphism are
  direct sums of $2^{m-1}$ terms on each of which the lemma for $m=1$
  can be applied.
\end{proof}

\begin{lem}\label{lem:equG}
  The following equalities hold in $\gr$:
  \begin{alignat*}{2}
    [x,a]=x^2,&\qquad [x,b]=x^2,\\
    [x,c]=x(1,x^{-1})x,&\qquad [x,d]=(1,x),\\
    [x^2,a]=x^4=((U,V)x^2,(V,U)x^2),&\qquad [x^2,b]=x^4,\\
    [x^2,c]=((U,V)x^2,(1,x)),&\qquad [x^2,d]=(1,(U,1)x^2),
  \end{alignat*}
  where $U=(1,x^{-1})x$ and $V=(x^{-1},1)x^{-1}$ are in $K$.
\end{lem}
\begin{proof}
  Direct computation; see also~\cite{rozhkov:habilitation}, where
  different notations are used.
\end{proof}

\begin{lem}\label{lem:techG}
  If $Q\gneqq N_{m+1}$ contains $g=(x,\dots,x)\in K_m$, then
  $[Q,\gr]\ge N_{m+1}$.
\end{lem}
\begin{proof}
  Let $b_m\in\{b,c,d\}$ be such that it acts like $b$ on
  $1^m\Sigma^*$. Then
  $$h=[g,b_m]=(1,\dots,1,[x,b])_m=(1,\dots,1,x^2)_m\in T_m.$$
  Conjugating $h$ by elements of $g$ yields all cyclic permutations of
  the above vector, so as $[G,\gr]$ is normal in $G$ it contains
  $T_m$. Likewise, let $d_m$ act like $d$ on $1^m\Sigma^*$. Then
  $$[g,d_m]=(1,\dots,1,[x,a],[x,d])_m=(1,\dots,1,x^2,(1,x))_m;$$
  using $T_m\le[Q,\gr]$, we obtain
  $(1,\dots,1,(1,x))_m=(1,\dots,1,x)_{m+1}\in[Q,\gr]$, so by the same
  conjugation argument $[Q,\gr]\ge K_{m+1}$.
\end{proof}


\newpage
\begin{thm} For all $m\ge1$ we have:
  \begin{enumerate}
  \item $\gamma_{2^m+1}(\gr)=N_m$.
  \item $\gamma_{2^m+1+r}(\gr)=N_{m+1}\alpha^{-1}(V_m^r)\beta^{-1}(V_{m-1}^r)$
    for $r=0,\dots,2^m$.
  \item $$\rank(\gamma_n(\gr)/\gamma_{n+1}(\gr))=\begin{cases}
      3 & \text{ if }n=1,\\
      2 & \text{ if }n=2,\\
      2 & \text{ if }n=2^m+1+r,\text{ with }0\le r<2^{m-1},\\
      1 & \text{ if }n=2^m+1+r,\text{ with }2^{m-1}\le r\le2^m.
    \end{cases}$$
  \end{enumerate}
\end{thm}

\begin{proof}
  First compute $\gamma_2(\gr)=\gr'=\langle[a,d],K\rangle$; it is of
  index $8$ in $\gr$, with quotient generated by $\{a,b,c\}$. Compute
  also $\gamma_3(\gr)=\langle x^2=[x,a],(1,x)=[x,d]\rangle^\gr=N_1$ of
  index $2$ in $\gamma_2(\gr)$, with quotient generated by
  $\{x^2,(1,x)\}$.  This gives the basis of an induction on $m\ge1$
  and $0<r\le2^m$.

  Assume that $\gamma_{2^m+1}(\gr)=N_m$. Note that the hypothesis of
  Lemma~\ref{lem:dec1} is satisfied; indeed $g_m$ can even be chosen
  among the conjugates of $b$, $c$ or $d$. Consider the sequence of
  quotients $Q_r=N_{m+1}\gamma_{2^m+1+r}(\gr)/N_{m+1}$ for $r\ge0$.
  Lemmata~\ref{lem:isoG} and~\ref{lem:dec1} tell us that
  $Q_r=\alpha^{-1}(V_m^r)\oplus\beta^{-1}(V_{m-1}^r)$; in particular
  $Q_r\ni(x,\dots,x)=\alpha^{-1}(v_m^{2^m-1})$ for all $r<2^m$, and
  then Lemma~\ref{lem:techG} tells us that $\gamma_{2^m+1+r}(\gr)\ge
  N_{m+1}$ for $r\le2^m$. When $r=2^m$ we have
  $\gamma_{2^{m+1}+1}(\gr)=N_{m+1}$ and the induction can continue.
\end{proof}
  
\begin{lem}\label{lem:squareG}
  For all $m\ge1$ and $r\in\{0,\dots,2^m-1\}$ we have:
  \begin{align*}
    (\alpha^{-1}V_m^r)^2&=\beta^{-1}(V_m^r)\le N_{m+1};\\
    (\beta^{-1}V_{m-1}^r)^2&=\beta^{-1}(V_m^{r+2^{m-1}})\le N_{m+1}.
  \end{align*}
\end{lem}
\begin{proof}
  Write $\alpha^{-1}(v_m^r)=(x^{i_1},\dots,x^{i_{2^m}})$ or
  $\beta^{-1}(v_m^r)=(x^{2i_1},\dots,x^{2i_{2^m}})$ for some
  $i_*\in\{0,1\}$. Then these
  claims follow immediately, using Lemma~\ref{lem:equG}, from
  \begin{align*}
    (\alpha^{-1}v_m^r)^2&=(x^{i_1},\dots,x^{i_{2^m}})^2=(x^{2i_1},\dots,x^{2i_{2^m}})=\beta^{-1}(v_m^r),\\
    (\beta^{-1}v_{m-1}^r)^2&=(x^{2i_1},\dots,x^{2i_{2^{m-1}}})^2=(x^{4i_1},\dots,x^{4i_{2^{m-1}}})\\
    &\equiv(x^{2i_1},x^{2i_1},\dots,x^{2i_{2^{m-1}}},x^{2i_{2^{m-1}}})=\beta^{-1}(v_m^{r+2^{m-1}})\mod N_{m+1}.
  \end{align*}
\end{proof}

\begin{thm}\label{thm:indG}
  For all $m\ge1$ we have:
  \begin{enumerate}
  \item $\gr_{2^m+1}=N_m$.
  \item $$\gr_{2^m+1+r}=\begin{cases}
      N_{m+1}\alpha^{-1}(V_m^r)\beta^{-1}(V_{m-1}^{r/2}) & \text{ if
        }0\le r\le2^m\text{ is even},\\
      N_{m+1}\alpha^{-1}(V_m^r)\beta^{-1}(V_{m-1}^{(r-1)/2}) & \text{ if
        }0\le r\le2^m\text{ is odd}.
    \end{cases}$$
  \item $$\rank(\gr_i/\gr_{i+1})=\begin{cases}
      3 & \text{ if }i=1,\\
      2 & \text{ if }i>1\text{ is even},\\
      1 & \text{ if }i>1\text{ is odd}.
    \end{cases}$$
  \end{enumerate}
\end{thm}
\begin{proof}
  First compute $\gr_2=\gamma_2(\gr)$ and $\gr_3=\gamma_3(\gr)=N_1$.
  This gives the basis of an induction on $m\ge1$ and $0\le r\le2^m$.
  Assume $\gr_{2^m+1}=N_m$.  Consider the sequence of quotients
  $Q_{m,r}=N_{m+1}\gr_{2^m+1+r}/N_{m+1}$ for $r\ge0$. We have
  $Q_{m,r}=[\gr,Q_{m,r-1}]Q_{m-1,\lfloor r/2\rfloor}^2$
  by~(\ref{eq:recj}).  Lemmata~\ref{lem:isoG}, \ref{lem:squareG}
  and~\ref{lem:dec1} tell us that
  $Q_r=\alpha^{-1}(V_m^r)\oplus\beta^{-1}(V_{m-1}^{\lfloor
    r/2\rfloor})$; in particular
  $Q_r\ni(x,\dots,x)=\alpha^{-1}(v_m^{2^m-1})$ for all $r<2^m$, and
  then Lemma~\ref{lem:techG} tells us that $\gr_{2^m+1+r}\ge N_{m+1}$
  for $r\le2^m$. When $r=2^m$ we have $\gr_{2^{m+1}+1}=N_{m+1}$ and the
  induction can continue.
\end{proof}

\subsection{Cayley graphs of Lie algebras}\label{subs:cayG}
We introduce the notion of \emph{Cayley graph} for graded Lie
algebras. Let $L=\bigoplus_{n=1}^\infty L_n$ be a graded Lie algebra
generated by a finite set $S$ of degree one elements. Fix a basis
$(\ell_{n,1},\dots,\ell_{n,\dim L_n})$ of $L_n$ for every $n$, and
give each $L_n$ an orthogonal scalar product
$\langle\ell_{n,i}|\ell_{n,j}\rangle=\delta_{i,j}$. The \emdef{Cayley
  graph} of $L$ is defined as follows: its vertices are the
$(i,j)\in\N^2$ with $i\ge1$ and $1\le j\le\dim L_n$. For every $s\in
S$ and $i,j,k\in\N$ there is an edge from $(i,j)$ to $(i+1,k)$
labeled by $s$ and with weight $\langle[\ell_{i,j},s]|\ell_{i+1,k}\rangle$.
By convention edges of weight $0$ are not represented. Additionally,
if $L$ is a $p$-algebra, there is an unlabeled edge of length
$(p-1)i$ from $(i,j)$ to $(pi,k)$ with weight $\langle
\ell_{i,j}^p|\ell_{pi,k}\rangle$.

Clearly, the Cayley graph of a Lie algebra $L$ determines the
structure of $L$. It is a connected graph, because $S$ is a generating set.
The geometric growth of the graph is the same as the growth of the
algebra.

As a simple example, consider the quaternion group $Q=\{\pm1,\pm i,\pm
j,\pm k\}$ generated by $\{i,j\}$, and its dimension series
$Q_1=Q$, $Q_2=\{\pm1\}$ and $Q_3=1$ over the field $\F_2$.
Then the Cayley graph of $\Lie(Q)$ is
$$
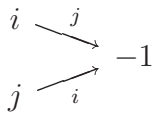
\begin{diagram}
  \node{i}\arrow{se,t}{j}\\
  \node[2]{-1}\\
  \node{j}\arrow{ne,b}{i}
\end{diagram}$$

We now describe the Cayley graphs of $L$ and $\Lie_{\F_2}$ associated
respectively to the lower central and dimensional series of $\gr$. Fix
$S=\{a,b,c,d\}$ as a generating set for $\gr$, and extend it to
$\overline
S=\{a,b,c,d,\{\begin{smallmatrix}b\\c\end{smallmatrix}\},\{\begin{smallmatrix}b\\d\end{smallmatrix}\},\{\begin{smallmatrix}c\\d\end{smallmatrix}\}\}$.
Define the transformation $\sigma$ on ${\overline S}^*$ by
$$\sigma(a)=a\{\begin{smallmatrix}b\\ c\end{smallmatrix}\}a,\quad\sigma(b)=d,\sigma(c)=b,\sigma(d)=b,$$
naturally extended to subsets. (For any fixed $g\in G$, one may obtain
all elements $h\in\stab_\gr(1)$ with $\phi(h)=(g,*)$ by computing
$\sigma(g)$ and making all possible choices of a letter from the
braced symbols. This explains the definition of $\overline S$.)

\begin{thm}\label{thm:lieG}
  The Cayley graph of $L(\gr)$ is as follows:

  $$
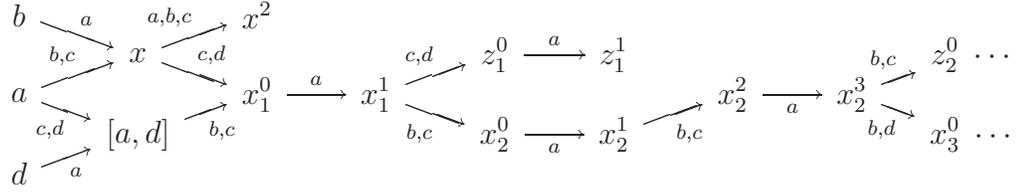
\begin{diagram}
    \node{b}\arrow{se,t}{a} \node[2]{x^2}\\
    \node[2]{x}\arrow{ne,t}{a,b,c}\arrow{se,t}{c,d}
    \node[3]{z_1^0}\arrow{e,t}{a} \node{z_1^1} \node[3]{z_2^0\;\cdots}
\\ 
    \node{a}\arrow{ne,t}{b,c}\arrow{se,b}{c,d}
    \node[2]{x_1^0}\arrow{e,t}{a}
    \node{x_1^1}\arrow{se,b}{b,c}\arrow{ne,t}{c,d}
    \node[3]{x_2^2}\arrow{e,b}{a}
    \node{x_2^3}\arrow{se,b}{b,d}\arrow{ne,t}{b,c}\\
    \node[2]{[a,d]}\arrow{ne,b}{b,c} \node[3]{x_2^0}\arrow{e,b}{a}
    \node{x_2^1}\arrow{ne,b}{b,c} \node[3]{x_3^0\;\cdots}
\\
    \node{d}\arrow{ne,b}{a}
  \end{diagram}$$
  where $x_m^r=\alpha^{-1}(v_m^r)$ and $z_m^r=\beta^{-1}(v_m^r)$. The
  edge \\
$(x_m^{2^m-1},x_{m+1}^0)$ is labelled by
  $\sigma^m\{\begin{smallmatrix}c\\d\end{smallmatrix}\}$, the edge
  $(x_m^{2^m-1},z_m^0)$ is labelled by
  $\sigma^m\{\begin{smallmatrix}b\\d\end{smallmatrix}\}$, and the
  paths from $x_m^0$ to $x_m^{2^m-1}$ and from $z_m^0$ to
  $z_m^{2^m-1}$ are labelled by $\sigma^{m-1}(a)$.

  The Cayley graph of $\Lie_{\F_2}(\gr)$ is as follows:
  $$
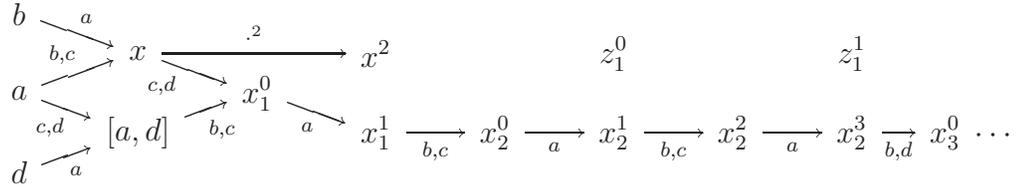
\begin{diagram}
    \node{b}\arrow{se,t}{a}\\
    \node[2]{x}\arrow{se,b,1}{c,d}\arrow[2]{e,t}{\cdot^2} \node[2]{x^2} \node[2]{z_1^0} \node[2]{z_1^1} 
\\
    \node{a}\arrow{ne,t}{b,c}\arrow{se,b}{c,d}
    \node[2]{x_1^0}\arrow{se,b}{a}\\
    \node[2]{[a,d]}\arrow{ne,b}{b,c}
    \node[2]{x_1^1}\arrow{e,b}{b,c}
    \node{x_2^0}\arrow{e,b}{a}
    \node{x_2^1}\arrow{e,b}{b,c}
    \node{x_2^2}\arrow{e,b}{a}
    \node{x_2^3}\arrow{e,b}{b,d}
    \node{x_3^0\;\cdots}
\\
    \node{d}\arrow{ne,b}{a}
  \end{diagram}$$
  with the same rule for labellings as for $L(\gr)$; and power maps
  from $x_m^r$ to $z_m^r$.
\end{thm}
Note that as we are in characteristic $2$ the non-zero weights can
only be $1$ and thus are not indicated.

\section{The Group $\grt$}\label{sec:S}
We describe here the lower central and dimension series for a
group $\grt$ containing the previous section's group $\gr$ as a
subgroup. More details about $\grt$ can be found
in~\cite{bartholdi-g:parabolic}.

As in Section~\ref{sec:G} set $\Sigma=\F_2$, and define automorphisms
$\tilde b$, $\tilde c$ and $\tilde d$ of $\Sigma^*$ by
\begin{alignat*}{2}
  \tilde b(0x\sigma)&=0\overline x\sigma,&\qquad\tilde b(1\sigma)&=1\tilde c(\sigma),\\
  \tilde c(0\sigma)&=0\sigma,&\qquad\tilde c(1\sigma)&=1\tilde d(\sigma),\\
  \tilde d(0\sigma)&=0\sigma,&\qquad\tilde d(1\sigma)&=1\tilde b(\sigma).
\end{alignat*}
Note that all generators are of order $2$ and $\{\tilde b,\tilde
c,\tilde d\}$ generate the elementary abelian group $2^3$. Set
$\grt=\langle a,\tilde b,\tilde c,\tilde d\rangle$. Clearly,
$\gr=\langle a,b=\tilde b\tilde c,c=\tilde c\tilde d,d=\tilde d\tilde
b\rangle$ is a subgroup of $\grt$. Its index is infinite, because
$\gr$ is a torsion group while $w=a\tilde b\tilde c\tilde d$ has
infinite order, because $w^2=(w^a,w)$.  Set $x=[a,\tilde b]$,
$y=[a,\tilde d]$, and
$$\tilde K = \langle x,y\rangle^\grt.$$
Then $\tilde K$ is a subgroup of finite index (actually index $32$) in $\gr$,
and contains $\tilde K\times\tilde K$ as a 
subgroup of finite index (actually index $8$). Set also $\tilde T=\langle
x^2\rangle^\gr=\tilde K^2$, and for any $Q\le\tilde K$ define
$Q_m=Q\times\dots\times Q$ ($2^m$ copies).  For $m\ge1$ set
$\tilde N_m=\tilde K_m\cdot\tilde T_{m-1}$.

For $m\ge2$, we have $\rist_\grt(m)=\tilde K_{m-2}$, so $\grt$ is a
branch group.

\begin{lem}\label{lem:isoS} The mapping
  $$\alpha\oplus\beta\oplus\gamma:\tilde N_m/\tilde
  N_{m+1}\longrightarrow V_m\oplus V_m\oplus V_{m-1}$$
  is an isomorphism for all $m$, where the $V_m$ are the modules
  defined in Subsection~\ref{subs:modules}, $\alpha$ maps
  $(1,\dots,x,\dots,1)\in\tilde K_m$ to the monomial in $V_m$
  corresponding to the vertex in $x$'s position, and $\beta$ maps
  $(1,\dots,y,\dots,1)\in\tilde K_m$ to the
  corresponding vertex in $V_m$, and $\gamma$ maps
  $(1,\dots,x^2,\dots,1)\in\tilde T_{m-1}$ to the corresponding
  monomial in $V_{m-1}$.
\end{lem}
\begin{proof}
  We first suppose $m=1$. Then 
$$\tilde N_1/\tilde N_2=\langle
  x^2,(1,x),(x,1),(1,y),(y,1)\rangle/\tilde N_2;$$ it is easy to check
  that $x^4=1$, so all generators of $\tilde N_1/\tilde N_2$ are of
  order $2$. Further, all commutators of generators belong to $\tilde
  K_2$, so the quotient $\tilde N_1/\tilde N_2$ is the elementary
  abelian group $2^5$, and $\alpha\oplus\beta\oplus\gamma$ is an
  isomorphism in that case.

  For $m>1$ it suffices to note that both sides of the isomorphism are
  direct sums of $2^{m-1}$ terms on each of which the lemma for $m=1$
  can be applied.
\end{proof}

\begin{lem}\label{lem:equS}
  The following equalities hold in $\grt$:
  \begin{alignat*}{2}
    [x,a]=x^2,&\qquad [x,\tilde b]=x^2\\
    [x,\tilde c]=(1,y),&\qquad [x,\tilde d]=(1,x),\\
    [x^2,a]=1,&\qquad [x^2,\tilde b]=1,\\
    [x^2,\tilde c]=1,&\qquad [x^2,\tilde d]=(1,x(x,1)x),\\
    [y,a]=1,&\qquad [y,\tilde b]=(x^{-1},1),\\
    [y,\tilde c]=1,&\qquad [y,\tilde d]=1.
  \end{alignat*}
\end{lem}
\begin{proof}
  Direct computation.
\end{proof}

\begin{lem}\label{lem:techS}
  If $Q\gneqq\tilde N_{m+1}$ contains $g=(x,\dots,x)\in\tilde K_m$, then
  $[Q,\grt]\ge\tilde N_{m+1}$.
\end{lem}
\begin{proof}
  Let $\tilde b_m\in\{\tilde b,\tilde c,\tilde d\}$ be such that it
  acts like $\tilde b$ on $1^m\Sigma^*$. Then
  $$[g,\tilde b_m]=(1,\dots,1,[x,\tilde b])_m=(1,\dots,1,x^2)_m\in\tilde T_m,$$
  so by a conjugation argument $[Q,\grt]\ge\tilde T_m$.
  Likewise, let $\tilde c_m$ and $\tilde d_m$ act like $c$ and $d$ on
  $1^m\Sigma^*$. Then
  \begin{align*}
    [g,\tilde c_m]&=(1,\dots,1,[x,a],[x,\tilde c])_m=(1,\dots,1,x^2,(1,y))_m,\\
    [g,\tilde d_m]&=(1,\dots,1,(1,x))_m.
  \end{align*}
  Using $\tilde T_m\le[Q,\grt]$, we obtain
  $(1,\dots,1,(1,y))_m=(1,\dots,1,y)_{m+1}\in[Q,\grt]$, so again by a
  conjugation argument $[Q,\grt]\ge\tilde K_{m+1}$.
\end{proof}

\begin{thm} For all $m\ge1$ we have:
  \begin{enumerate}
  \item $\gamma_{2^m+1}(\grt)=\tilde N_m$.
  \item $\gamma_{2^m+1+r}(\grt)=\tilde N_{m+1}\alpha^{-1}(V_m^r)\beta^{-1}(V_m^r)\gamma^{-1}(V_{m-1}^r)$
    for $r=0,\dots,2^m$.
  \item $$\rank(\gamma_n(\grt)/\gamma_{n+1}(\grt))=\begin{cases}
      4 & \text{ if }n=1,\\
      3 & \text{ if }n=2,\\
      3 & \text{ if }n=2^m+1+r,\text{ with }0\le r<2^{m-1},\\
      2 & \text{ if }n=2^m+1+r,\text{ with }2^{m-1}\le r\le2^m.
    \end{cases}$$
  \end{enumerate}
\end{thm}
\begin{proof}
  First compute $\gamma_2(\grt)=\grt'=\langle[a,\tilde
  c],\tilde K\rangle$, of index $16$ in $\grt$, and
  $\gamma_3(\grt)=\langle x^2,(1,x),(1,y)\rangle^\grt=\tilde
  N_1$, with $x^2=[x,a]$, $(1,x)=[x,\tilde d]$ and~$(1,y)=[x,\tilde c]$.
  This gives the basis of an induction on $m\ge1$ and $0\le r\le2^m$.
  
  Assume that $\gamma_{2^m+1}(\grt)=\tilde N_m$. Note that the
  hypothesis of Lemma~\ref{lem:dec1} is satisfied for $\grt$, as it
  holds for $\gr<\grt$. Consider the sequence of quotients $Q_r=\tilde
  N_{m+1}\gamma_{2^m+1+r}(\grt)/\tilde N_{m+1}$ for $r\ge0$.
  Lemmata~\ref{lem:isoS} and~\ref{lem:dec1} tell us that
  $Q_r=\alpha^{-1}(V_m^r)\oplus\beta^{-1}(V_m^r)\oplus\gamma^{-1}(V_{m-1}^r)$;
  in particular $Q_r\ni(x,\dots,x)=\alpha^{-1}(v_m^{2^m-1})$ for all
  $r<2^m$, and then Lemma~\ref{lem:techS} tells us that
  $\gamma_{2^m+1+r}(\grt)\ge\tilde N_{m+1}$ for $r\le2^m$. When
  $r=2^m$ we have $\gamma_{2^{m+1}+1}(\grt)=\tilde N_{m+1}$ and the
  induction can continue.
\end{proof}
  
\begin{lem}\label{lem:squareS}
  For all $m\ge1$ and $r\in\{0,\dots,2^m-1\}$ we have:
  \begin{align*}
    (\alpha^{-1}V_m^r)^2&=\gamma^{-1}(V_m^r)\le \tilde N_{m+1};\\
    (\beta^{-1}V_m^r)^2&=1\le\tilde N_{m+1};\\
    (\gamma^{-1}V_{m-1}^r)^2&=\gamma^{-1}(V_m^{r+2^{m-1}})\le\tilde N_{m+1}.
  \end{align*}
\end{lem}
\begin{proof}
  Write $\alpha^{-1}(v_m^r)=(x^{i_1},\dots,x^{i_{2^m}})$,
  $\beta^{-1}(v_m^r)=(y^{i_1},\dots,y^{i_{2^m}})$ or\\
  $\gamma^{-1}(v_m^r)=(x^{2i_1},\dots,x^{2i_{2^m}})$ for some
  $i_*\in\{0,1\}$. Then these claims follow immediately, using
  Lemma~\ref{lem:equS}, from
  \begin{align*}
    (\alpha^{-1}v_m^r)^2&=(x^{i_1},\dots,x^{i_{2^m}})^2=(x^{2i_1},\dots,x^{2i_{2^m}})=\gamma^{-1}(v_m^r),\\
    (\beta^{-1}v_m^r)^2&=(y^{i_1},\dots,y^{i_{2^m}})^2=(y^{2i_1},\dots,y^{2i_{2^m}})=(1,\dots,1),\\
    (\gamma^{-1}v_{m-1}^r)^2&=(x^{2i_1},\dots,x^{2i_{2^{m-1}}})^2=(x^{4i_1},\dots,x^{4i_{2^{m-1}}})\\
    &\equiv(x^{2i_1},x^{2i_1},\dots,x^{2i_{2^{m-1}}},x^{2i_{2^{m-1}}})=\gamma^{-1}(v_m^{r+2^{m-1}})\mod\tilde N_{m+1}.
  \end{align*}
\end{proof}

\begin{thm}\label{thm:indS}
  For all $m\ge1$ we have:
  \begin{enumerate}
  \item $\grt_{2^m+1}=\tilde N_m$.
  \item $$\grt_{2^m+1+r}=\begin{cases}
      \tilde N_{m+1}\alpha^{-1}(V_m^r)\beta^{-1}(V_m^r)\gamma^{-1}(V_{m-1}^{r/2}) & \text{ if }0\le r\le2^m\text{ is even},\\
      \tilde N_{m+1}\alpha^{-1}(V_m^r)\beta^{-1}(V_m^r)\gamma^{-1}(V_{m-1}^{(r-1)/2}) & \text{ if }0\le r\le2^m\text{ is odd}.
    \end{cases}$$
  \item $$\rank(\grt_i/\grt_{i+1})=\begin{cases}
      4 & \text{ if }i=1,\\
      3 & \text{ if }i>1\text{ is even},\\
      2 & \text{ if }i>1\text{ is odd}.
    \end{cases}$$
  \end{enumerate}
\end{thm}
\begin{proof}
  First compute $\grt_2=\gamma_2(\grt)$ and
  $\grt_3=\gamma_3(\grt)=\tilde N_1$. This gives the basis of an
  induction on $m\ge1$ and $0\le r\le2^m$. Assume $\gr_{2^m+1}=\tilde
  N_m$.  Consider the sequence of quotients $Q_{m,r}=\tilde
  N_{m+1}\grt_{2^m+1+r}/\tilde N_{m+1}$ for $r\ge0$. We have
  $Q_{m,r}=[\grt,Q_{m,r-1}]Q_{m-1,\lfloor r/2\rfloor}^2$
  by~(\ref{eq:recj}).  Lemmata~\ref{lem:isoS}, \ref{lem:squareS}
  and~\ref{lem:dec1} tell us that
  $Q_r=\alpha^{-1}(V_m^r)\oplus\beta^{-1}(V_m^r)\oplus\gamma^{-1}(V_{m-1}^{\lfloor r/2\rfloor})$; in particular
  $Q_r\ni(x,\dots,x)=\alpha^{-1}(v_m^{2^m-1})$ for all $r<2^m$, and
  then Lemma~\ref{lem:techS} tells us that $\gr_{2^m+1+r}\ge\tilde
  N_{m+1}$ for $r\le2^m$. When $r=2^m$ we have $\gr_{2^{m+1}+1}=\tilde
  N_{m+1}$ and the induction can continue.
\end{proof}

\subsection{The Lie Algebra Structures}
We describe here the Cayley graphs of $L$ and $\Lie_\Fp$ associated
respectively to the lower central and dimension series of $\grt$.
Consider $\tilde S=\{a,\tilde b,\tilde c,\tilde d\}$ and define the
transformation $\tilde\sigma$ on ${\tilde S}^*$ by
$$\tilde\sigma(a)=a\tilde ba,\quad\tilde\sigma(\tilde b)=\tilde d,\quad\tilde\sigma(\tilde c)=\tilde b,\quad\tilde\sigma(\tilde d)=\tilde b.$$

\begin{thm}\label{thm:lieS}
  The Cayley graph of $L(\grt)$ is as follows:
  $$
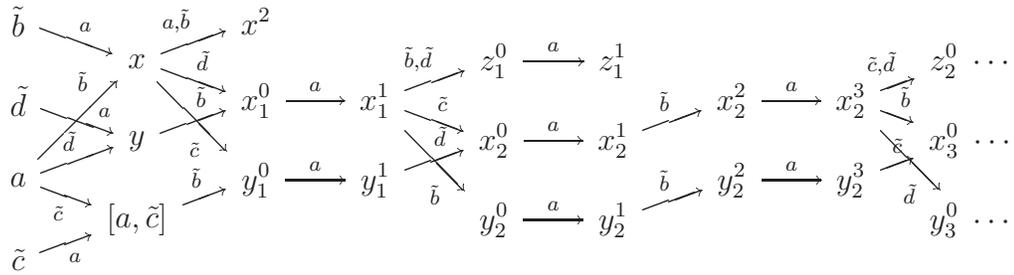
\begin{diagram}
    \node{\tilde b}\arrow{se,t}{a} \node[2]{x^2}\\
    \node[2]{x}\arrow{ne,t}{a,\tilde b}\arrow{se,t}{\tilde d}\arrow{ssse,b,3}{\tilde c} \node[3]{z_1^0}\arrow{e,t}{a} \node{z_1^1} \node[3]{z_2^0\;\cdots}
\\
    \node{\tilde d}\arrow{se,t,3}{a} \node[2]{x_1^0}\arrow{e,t}{a} \node{x_1^1}\arrow{ne,t}{\tilde b,\tilde d}\arrow{se,t}{\tilde c}\arrow{ssse,b,3}{\tilde b} \node[3]{x_2^2}\arrow{e,t}{a} \node{x_2^3}\arrow{ne,t}{\tilde c,\tilde d}\arrow{se,t}{\tilde b}\arrow{ssse,b,3}{\tilde d}\\
    \node[2]{y}\arrow{ne,t,3}{\tilde b} \node[3]{x_2^0}\arrow{e,t}{a} \node{x_2^1}\arrow{ne,t}{\tilde b} \node[3]{x_3^0\;\cdots}
\\
    \node{a}\arrow{ne,t}{\tilde d}\arrow{se,b}{\tilde c}\arrow{nnne,t,3}{\tilde b} \node[2]{y_1^0}\arrow{e,t}{a} \node{y_1^1}\arrow{ne,t,3}{\tilde d} \node[3]{y_2^2}\arrow{e,t}{a} \node{y_2^3}\arrow{ne,t,3}{\tilde c}\\
    \node[2]{[a,\tilde c]}\arrow{ne,t}{\tilde b} \node[3]{y_2^0}\arrow{e,t}{a} \node{y_2^1}\arrow{ne,t}{\tilde b} \node[3]{y_3^0\;\cdots}
\\
    \node{\tilde c}\arrow{ne,b}{a}
  \end{diagram}$$
  where $x_m^r=\alpha^{-1}(v_m^r)$, $y_m^r=\beta^{-1}(v_m^r)$ and
  $z_m^r=\gamma^{-1}(v_m^r)$. The edge \\
$(x_m^{2^m-1},x_{m+1}^0)$ is
  labelled by $\tilde\sigma^m(\tilde d)$, the edges
  $(x_m^{2^m-1},y_{m+1}^0)$ and $(x_m^{2^m-1},z_m^0)$ are labelled by
  $\tilde\sigma^m(\tilde b)$, the edges $(x_m^{2^m-1},z_m^0)$ and
  $(y_m^{2^m-1},x_{m+1}^0)$ are labelled by $\tilde\sigma^m(\tilde
  c)$, and the paths from $x_m^0$ to $x_m^{2^m-1}$, from $y_m^0$ to
  $y_m^{2^m-1}$ and from $z_m^0$ to $z_m^{2^m-1}$ are labelled by
  $\tilde\sigma^{m-1}(a)$.

  The Cayley graph of $\Lie_{\F_2}(\grt)$ is as follows:
  $$
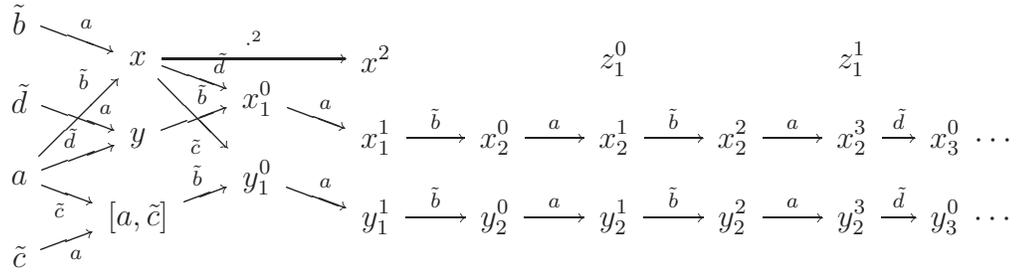
\begin{diagram}
    \node{\tilde b}\arrow{se,t}{a}\\
    \node[2]{x}\arrow{se,t,3}{\tilde d}\arrow{ssse,b,3}{\tilde c}\arrow[2]{e,t}{\cdot^2} \node[2]{x^2} \node[2]{z_1^0} \node[2]{z_1^1} 
\\
    \node{\tilde d}\arrow{se,t,3}{a} \node[2]{x_1^0}\arrow{se,t}{a}\\
    \node[2]{y}\arrow{ne,t,3}{\tilde b} \node[2]{x_1^1}\arrow{e,t}{\tilde b} \node{x_2^0}\arrow{e,t}{a} \node{x_2^1}\arrow{e,t}{\tilde b} \node{x_2^2}\arrow{e,t}{a} \node{x_2^3}\arrow{e,t}{\tilde d} \node{x_3^0\;\cdots}
\\
    \node{a}\arrow{ne,t}{\tilde d}\arrow{se,b}{\tilde c}\arrow{nnne,t,3}{\tilde b} \node[2]{y_1^0}\arrow{se,t}{a}\\
    \node[2]{[a,\tilde c]}\arrow{ne,t}{\tilde b} \node[2]{y_1^1}\arrow{e,t}{\tilde b} \node{y_2^0}\arrow{e,t}{a} \node{y_2^1}\arrow{e,t}{\tilde b} \node{y_2^2}\arrow{e,t}{a} \node{y_2^3}\arrow{e,t}{\tilde d} \node{y_3^0\;\cdots}
\\
    \node{\tilde c}\arrow{ne,b}{a}
  \end{diagram}$$
  with the same labellings as for $L(\grt)$; and power maps from $x_m^r$ to
  $z_m^r$.
\end{thm}

\section{Other Fractal Groups}
The technique involved in the proof of the results of the last three
sections show that for a group $G$ acting on a tree $\Sigma^*$ by powers
of the cyclic permutation $\epsilon=(0,1,\dots,p-1)$ at each vertex,
$G$ has finite width when the following conditions are satisfied:
\begin{enumerate}
\item the corresponding action on a sequence $\{V_n\}_{n=0}^\infty$ of
  $G$-modules as defined in Subsection~\ref{subs:modules} has the
  \emdef{bounded corank property}, i.e.\ there is a constant $C$ such that
  $$\dim V_n^r/[G,V_n^r]\le C$$
  for all $n\ge0$ and $0\le r\le p^n-1$.
\item There is a descending sequence $\{N_m\}_{m=1}^\infty$ of normal
  subgroups of $G$ satisfying the condition that for all $m$ the
  quotients $N_m/N_{m+1}$ are isomorphic to some direct sum
  $\bigoplus_{i=1}^KV_{m+\delta_i}$ for fixed $K$ and $\delta_i$.
\end{enumerate}

Let us mention that the $p$-groups $G_\omega$, for arbitrary $p\ge2$
and $\omega\in\{0,\dots,p\}^\N$ constructed
in~\cite{grigorchuk:gdegree,grigorchuk:pgps} all satisfy
Condition~1. Also, the group $\langle a,t\rangle<\aut(\Sigma_p^*)$,
$p\ge3$, where $a$ acts as $\epsilon$ on the root vertex and trivially
elsewhere and $t$ is defined recursively by $t=(a,1,\dots,1,t)$,
satisfies Condition~1. We believe that this last group also satifies
Condition~2, as do all $G_\omega$ for periodic sequences
$\omega$. Note that $\gr$ is a particular case of $G_\omega$ when
$p=2$ and $\omega=012012\dots$. Therefore they all `should' have
finite width.

Meanwhile, the Gupta-Sidki groups constructed in~\cite{gupta-s:3group} do
not satisfy Condition~1. As was proved recently by the first author,
the growth of the Lie algebra $\Lie_\Fp(G)$ coincides with the
spherical growth of the Schreier graph of $G$ relatively to
$\stab_G(e)$, where $e$ is an infinite geodesic path in the tree
$\Sigma^*$. For our groups $\gr$ and $\grt$ the spherical growth is
bounded and this is why these groups have bounded width. For the
Gupta-Sidki groups, the spherical growth of the Schreier graph is
unbounded (it grows approximately as $\sqrt n$), and therefore these
groups do not have the finite width property. It also follows from
these considerations that their growth is at least
$e^{n^{1-1/(1/2+2)}}=e^{n^{3/5}}$.

\section{Profinite Groups of Finite Width}
Finally we wish to explain how our results in the previous
sections lead to counterexamples to Conjecture~\ref{conj:zelmanov}
stated in the introduction. Let $\widehat G$ be the profinite
completion of $G=\gr$ or $\grt$.
\begin{thm}
  The group $\widehat G$ is a just-infinite pro-$2$-group of finite
  width which does not belong to the list of
  Conjecture~\ref{conj:zelmanov} (which consists of solvable groups,
  $p$-adic analytic groups, and groups commensurable to positive parts
  of loop groups or to the Nottingham group).
\end{thm}

Its proof relies on the following notion:
\begin{defn}
  Let $G<\aut(\Sigma^*)$ be a group acting on a rooted tree. $G$ has
  the \emdef{congruence subgroup property} if for any finite-index
  subgroup $H$ of $G$ there is an $n$ such that $\stab_G(n)<H<G$.
\end{defn}
\begin{proof}
  $G$ has the congruence property. This is well known for $\gr$ (see
  for instance~\cite{grigorchuk:jibg}); while for $\grt$ the subgroup
  $\tilde K$ contains $\stab_\grt(4)$ and enjoys the property that every
  subgroup of finite index in $\grt$ contains $\tilde K_m=\tilde
  K\times\dots\times\tilde K$ for some $m$;
  see~\cite{bartholdi-g:parabolic}.

  The profinite completion of $G$ with respect to its subgroups
  $\stab_G(n)$ is therefore a pro-$2$-group and coincides with the
  closure of $G$ in $\aut(\{0,1\}^*)$. The closure of a branch group
  is again a branch group, as is observed in~\cite{grigorchuk:jibg}.
  
  The criterion of just-infiniteness for profinite branch groups is
  the same as the one for discrete branch groups given
  in~\cite{bartholdi-g:parabolic}; it is that $\overline K/\overline
  K'$ (respectively $\overline{\tilde K}/\overline{\tilde K}'$) be
  finite, where $\overline K$ and $\overline{\tilde K}$ are the
  closures of $K$ and $\tilde K$. Now $|\overline K/\overline
  K'|\le|K/K'|<\infty$, the last inequality following from a
  computation in~\cite{bartholdi-g:parabolic}. The same inequalities
  hold for $\grt$, and this proves the just-infiniteness of $\widehat
  G$.
  
  The group $\widehat G$ has finite width for both versions of
  Definition~\ref{defn:fw}. This is clear for $D$-width, because the
  discrete and pro-$p$ Lie algebras $\Lie(G)$ and $\Lie(\widehat G)$
  are isomorphic. The finiteness of $C$-width follows from the
  inequalities
  $$\left|\gamma_n(\widehat G)/\gamma_{n+1}(\widehat G)\right|\le\left|\gamma_n(G)/\gamma_{n+1}(G)\right|<\infty,$$
  which again are consequences of the congruence property of $G$.
  
  Finally, $\widehat G$ does not belong to the list of groups given in
  Conjecture~\ref{conj:zelmanov}: it is neither solvable, because $G$
  isn't, nor $p$-adic analytic, by Lazard's
  criterion~\cite{lazard:padanal} (its Lie algebra $\Lie(\widehat
  G)=\Lie(G)$ would have a zero component in some dimension). The
  other groups in the list of Conjecture~\ref{conj:zelmanov} are
  \emdef{hereditarily just-infinite} groups, that is, groups every
  open subgroup of which is
  just-infinite~\cite[page~5]{klass-lg-p:fw}. Profinite just-infinite
  branch groups are never hereditarily just-infinite, as is shown
  in~\cite{grigorchuk:jibg}.
\end{proof}

\begin{tabbing}
Section de Math\'ematiques \hspace*{170pt} \= Department of Ordinary Differential Equations \kill
Laurent Bartholdi \> Rostislav Grigorchuk \\
Section de Math\a'ematiques \> Steklov Mathematical Institute\\
Universit\a'e de Gen\a`eve \> Gubkina Street 8\\
CP 240, 1211 Gen\a`eve 24 \> Moscow 117966\\
Switzerland \>  Russia\\
{\tt Laurent.Bartholdi@math.unige.ch} \> {\tt grigorch@mi.ras.ru}
\end{tabbing}

\bibliography{mrabbrev,people,math,grigorchuk,bartholdi}
\end{document}